\documentclass[12pt]{article}
\usepackage{latexsym,amsmath,amssymb,amsfonts,amsthm}
\usepackage[dvips]{graphicx}
\usepackage{epic,eepic}

\theoremstyle{plain}
\newtheorem{theorem}{Theorem }[section]
\newtheorem{proposition}[theorem]{Proposition}
\newtheorem{lemma}[theorem]{Lemma}
\newtheorem{corollary}[theorem]{Corollary}
\newtheorem{maintheorem}{Theorem}
\newtheorem{claim}[theorem]{Claim}

\theoremstyle{definition}
\theoremstyle{remark}
\newtheorem{remark}[theorem]{Remark}

\newtheorem{definition}[theorem]{Definition}

\newcommand{\field}[1]{\mathbb{#1}}
\newcommand{\real}{\field{R}}

\newcommand{\integer}{\field{Z}}

\newcommand{\Circ}{\field{S}}

\newcommand{\al} {\alpha}       
\newcommand{\be} {\beta}        
\newcommand{\ga} {\gamma}    
       
\newcommand{\ep} {\epsilon}

\newcommand{\ka} {\kappa}

\newcommand{\SC}{{\cal C}}

\newcommand{\SI}{{\cal I}}

\newcommand{\SK}{{\cal K}}

\newcommand{\SP}{{\cal P}}

\newcommand{\p}{\partial}

\begin{document}

\title{\Large{\textbf{ Maximal entropy measures for Viana maps }}}

\author{Alexander Arbieto, Carlos Matheus and Samuel Senti}

\date{May 30, 2005}

\maketitle

\begin{abstract}In this note we construct measures of maximal entropy for a 
certain class of maps with critical points called Viana maps. The main
ingredients of the proof are the non-uniform expansion features and the slow
recurrence (to the critical set) of generic points with respect to the natural 
candidates for attaining the topological entropy.
\end{abstract}

\section{Introduction}
Nowadays, it is well-known that a dynamical system can be understood from the
study of its invariant measures (this study is called \emph{ergodic
theory}). However, there are in general many invariant measures for a
given system, and it is necessary to make a selection of some ``interesting''
measures. Motivated by statistical physics, some candidates called
\emph{equilibrium states}, are the measures which satisfy a variational principle.

The theory of equilibrium states in the case of uniformly hyperbolic systems
(developed by Bowen, Ruelle, Sinai among others) is now a classical theory with
completely satisfactory results. But, in the non-uniform (higher dimensional) 
context only few results are known (see, for instance, Oliveira~\cite{O} for a
recent progress in this direction; for a random version of it,
see Arbieto, Matheus, Oliveira~\cite{AMO}).

The central theme of this note is to prove the existence of maximal entropy
measures for a class of non-uniformly expanding maps with critical set 
introduced by Viana~\cite{V}. Some features of these maps (which we would like 
to call \emph{Viana maps}) include the existence of a unique absolutely 
continuous invariant measure $\mu_0$ with non-uniform expansion, slow 
recurrence to the critical set and positive Lyapounov exponents (see~\cite{V}). 
Moreover, $\mu_0$ is stochastically stable\footnote{Here, the stochastic 
stability is restricted to the random perturbations with the same critical set
$\SC$ and whose derivative at any point $x\notin\SC$ (i.e., in the complement
of the critical set) is equal to the derivative of the unperturbed system at
the same point $x$ (see Alves, Ara\'ujo~\cite{AA})}. 

On the other hand, an obstruction to construct general equilibrium states for 
Viana maps (among other non-uniformly expanding maps) is the \emph{a priori} 
dependence on the Lebesgue measure as reference in the concepts of 
non-uniform expansion and slow recurrence to the critical set. To overcome 
this difficult we apply the basic strategy of Oliveira~\cite{O}: 

\begin{itemize}
\item We select a certain set $\SK$ of 
invariant measures which are natural candidates to realize the 
variational principle (for nearly constant potentials); 
\item We prove that 
any measure outside $\SK$ can not be an equilibrium measure and any measure 
inside $\SK$ is \emph{expanding} (see the definition~\ref{d.expanding} below);
\item We show that any expanding measure admits generating partitions. This
leads us to the semicontinuity of the entropy, and hence, by a standard
argument, to the existence of equilibrium states for nearly constant
potentials.
\end{itemize}

However, it is worth pointing out an extra difficulty of the case of Viana
maps: the presence of critical points is an obstacle to obtaining infinitely
many times of uniform expansion directly from the non-uniform expansion
condition (indeed, some slow recurrence to the critical set should be
proven).\footnote{In the context of Oliveira~\cite{O}, this problem does not
occur since the transformations are local diffeomorphisms, and so infintely
many hyperbolic times can be obtained directly from non-uniform expansion.}
Since the abundance of hyperbolic times is fundamental to prove that
expanding measures admit generating partitions (see lemma~\ref{l.partition}), 
we need to solve this problem.

Fortunately, the slow recurrence to the critical set can be obtained from the
integrability of the distance function to the critical set $\SC$. Because 
the Viana maps \emph{behave like a power of the distance} to $\SC$, it is
reasonable that this integrability problem is related to the regularity of the 
Lyapounov exponents (i.e., the Lyapounov exponents are bounded away from
$-\infty$). Now, the regularity of the Lyapounov exponents of measures with a
non-trivial chance to attain the supremum of the variational principle (i.e.,
whose entropy is not $-\infty$)
is an easy consequence of Ruelle's inequality. Thus, the presence of critical
points is not a great trouble and we still can use the basic strategy of
Oliveira~\cite{O}.

Now we are going to state our main result. In order to do this, let us recall
some definitions and notations.

In general, given a continuous map $f:M\to M$ of a compact metric space $M$ and
a continuous function $\phi:M\to\real$, we call an 
$f$-invariant Borel probability measure
$\mu$ an \emph{equilibrium state} of $(f,\phi)$ if $\mu$ realizes the
variational principle 
$$h_{\mu}(f)+\int\phi d\mu = \sup\limits_{\eta\in\SI} 
\left(h_{\eta}(f)+\int\phi d\eta \right),$$
where $\SI$ is the set of $f$-invariant Borel probabilities.

Now we are in position to state our main result:

\begin{maintheorem}\label{t.A}Viana maps admit equilibrium states for any
nearly constant potential $\phi$ (in the sense of the 
condition~(\ref{e.potential}) below). In particular, Viana maps have maximal
entropy measure (since they are equilibrium states for the constant potential
$\phi=0$). Moreover, any equilibrium state of the Viana maps (associated to
nearly constant potentials) are hyperbolic measures with all Lyapounov exponents
greater than some $\lambda_0>0$.   
\end{maintheorem}

To close the introduction, let us comment about the organization of the paper: 
section~\ref{s.Viana} contains the definition of Viana maps, some of its
properties and a precise statement of what does ``$\phi$ is nearly
constant'' means. In section~\ref{s.proof} we present the proof of the
theorem~\ref{t.A}. Finally, in section~\ref{s.remarks}, we point out some
generalizations and questions concernig the theorem~\ref{t.A}.

\section{Viana maps}\label{s.Viana}

Let $a_0\in (1,2)$ such that $x=0$ is pre-periodic for the quadratic map
$h(x)=a_0-x^2$. Denote by $\Circ^1=\real / \integer$ and $b:\Circ^1\to\real$ a
Morse function, e.g., $b(\theta)=\sin (2\pi\theta)$. Fix some $\al>0$
sufficiently small and put 
$$
\widetilde{f} (\theta,x) := (g (\theta), a(\theta)-x^2) := (g(\theta), 
Q(\theta,x)),
$$
where $g$ is the expanding map $g(\theta) = d\theta$ of $\Circ^1$, 
for some integer $d\geq 16$ and $a(\theta) = a_0+\al b(\theta)$.
\footnote{The assumption $d\geq 16$ doesn't play crucial role in our results
and it is present here only for sake of simplicity. Indeed, the results of 
Buzzi, Sester and Tsujii~\cite{BST} can be applied to replace $d\geq 16$ by 
$d\geq 2$, at least when we consider the $C^{\infty}$-topology. See the
section~\ref{s.remarks} for details.}   
Since $a_0 < 2$, it is easy to check that there is $I\subset (-2,2)$ such that 
the closure of $\widetilde{f} (\Circ^1\times I)$ is contaned in the interior
of $\Circ^1\times I$. Indeed, $h(x)=a_0-x^2$ has a unique fixed point
$x_0=\frac{1}{2}(-1-\sqrt{1+4a_0})<0$ which is a repeller.
In particular, if we take $\beta>0$ slightly smaller than $-x_0$, 
the interval $I= [ -\beta,\beta ]$ satisfies $h(I)\subset\text{int}(I)$ and
$|h'|>1$ on $\real\setminus\text{int}(I)$.
Hence, if $\al$ is sufficienlty small we still have 
$\widetilde{f} (\Circ^1\times I)\subset\text{int}(\Circ^1\times I)$.
Notice also that, if 
$f$ is $C^0$-close to $\widetilde{f}$ then $f$ also has 
$\Circ^1\times I$ as a forward invariant region, i.e., any $f$ close to
$\widetilde{f}$ may be considered as a map 
$f:\Circ^1\times I\to\Circ^1\times I$. 

In what follows, we consider the parameters $d\geq 16$ and $a_0\in (1,2)$
fixed and $\al$ is sufficiently small depending on $d$ and $a_0$. Under these
conditions, Viana proved that

\begin{theorem}[Viana~\cite{V}]\label{t.exponent}
For $\al$ is sufficiently small, there exists a positive constant $c_0>0$ such
that any map $f$ sufficiently close to $\widetilde{f}$ in the $C^3$ topology
has both Lyapounov exponents greater than $c_0$ at Lebesgue almost every
point, that is, 
$$
\liminf\limits_{n\to\infty}\frac{1}{n}\log\|Df^n(\theta,x)\| > c_0, \quad
\text{ for Lebesgue a.e. } (\theta,x)\in\Circ^1\times I.
$$
\end{theorem} 

A very common fact in smooth ergodic theory is: ``the positivity of Lyapounov 
exponents acts as a strong 
evidence of the existence of absolutely continuous invariant measures''. Indeed,
this general principle works in several examples of maps with positive
Lyapounov exponents, including Viana maps:

\begin{theorem}[Alves~\cite{A}]\label{t.SRB}
Any map $f$ sufficiently $C^3$ close to $\widetilde{f}$ admits an unique
absolutely continuous invariant measure $\mu_0$.
\end{theorem} 

Next, let us study the tangent bundle dynamics of the Viana maps. Denote by
$f_0$ the product map $f_0 = g\times h$, i.e., $f_0:\Circ^1\times I
\to\Circ^1\times I$, $f_0(\theta ,x) = (d\theta, a_0 - x^2)$. Since 
$Df_0\cdot\frac{\p}{\p\theta} = d\cdot\frac{\p}{\p\theta}$,
$Df_0\cdot\frac{\p}{\p x} = -2x\frac{\p}{\p x}$, $d\geq 16 > 4$ and $I\subset
(-2,2)$, it is not hard to see that the splitting $T_{(\theta,x)}
(\Circ^1\times I) = E^u\oplus E^c$ is a dominated decomposition of the tangent
bundle of $\Circ^1\times I$ into two invariant (one-dimensional) subbundles
where $E^u$ is uniformly expanding and $E^c$ is dominated by $E^u$. In other
words, $f_0$ is a partially hyperbolic map of type $E^u\oplus E^c$. From the
theory of partial hyperbolicity we know that, if $\al$ is sufficiently small,
$\widetilde{f}$ (as any $C^1$ nearby map) also admits a partially hyperbolic
splitting of type $E^u\oplus E^c$ which varies continuously. As a consequence,
we get 

\begin{proposition}\label{p.conformal}Given $\ep>0$, the Lyapounov exponent 
$$\lambda^u (p,f):= 
\lim\limits_{n\to\pm\infty}\frac{1}{n}\log\|Df^n(p)|_{E^u}\|$$ 
associated to $E^u$ at every point $p$ verifies 
\begin{equation}\label{e.conformal}
\log (d-\ep)\leq\lambda^u (p,f)\leq \log (d+\ep),
\end{equation}
for any $f$ close 
to $\widetilde{f}$ (if $\al$ is sufficiently small). 
\end{proposition}

\begin{proof}Since $E^u$ is a continuous function of $f$ and the derivative of 
$f_0$ expands $\frac{\p}{\p\theta}$ by the constant factor $d$, it is clear that
the proposition holds for any $f$ sufficiently $C^1$-close to $f_0$, which is
certainly true in the context of Viana maps, provided $\al$ is small.  
\end{proof}  

Before ending this section, we apply the previous proposition to obtain a lower
bound on the entropy of the SRB measure $\mu_0$. Using Pesin's formula we know
that the entropy $h_{\mu_0}(f)$ of $\mu_0$ is the integrated sum of the positive
Lyapounov exponents of $\mu_0$. Hence, if we use theorem~\ref{t.exponent} and
proposition~\ref{p.conformal}, it follows that 
\begin{equation}\label{e.1}
h_{\mu_0}(f)\geq \log (d-\ep) + c_0.
\end{equation}

Once this notation is established, we are able to state precisely our condition
on the potential $\phi$. 

\begin{definition}
We say that the potential $\phi$ is \emph{nearly constant} if
\begin{equation}\label{e.potential}
\max\phi - \min\phi < 
\frac{1}{2}(c_0 - \frac{\log (d+\ep)}{\log (d-\ep)}).
\end{equation}
Note that the right-hand side of~(\ref{e.potential}) is positive for $\ep>0$
sufficiently small depending on $d$ and $c_0$, the lower bound on the Lyapounov
exponents of the Viana maps.
\end{definition}

After these preparation, we are ready to prove our main result.

\section{Proof of theorem~\ref{t.A}}\label{s.proof}

From now on, we consider Viana maps $f$, which are, by definition, all $C^3$
close maps to $\widetilde{f}$ for any small $\al$. 

We start with the control of the recurrence (to the critical set $\SC$) of 
generic points of invariant ergodic measures whose Lyapounov exponents are 
regular. The first step is to prove the following lemma:

\begin{lemma}\label{l.recurrence}Let $\eta$ be an ergodic measure such that
$\lambda^c(\eta) > -\infty$, where $\lambda^c$ is the Lyapounov exponent of
$\eta$ associated to $E^c$. Then, 
$$
\int |\log\,\text{dist}\,(p,\SC) | d\eta < \infty.
$$
\end{lemma}

\begin{proof}Since $\eta$ is ergodic, $\lambda^c (\eta)= 
\int\log\|Df|_{E^c}\| d\eta > -\infty$. On the other hand, the definition of 
$f_0$ and $\al$ sufficiently small implies 
$ \frac{1}{3}\|Df(p)|_{E^c}\|\leq dist(p,\SC)\leq 3\|Df(p)|_{E^c}\| $.
These two facts together finish the proof.
\end{proof}

A consequence of this lemma is the slow recurrence (to $\SC$) of generic points
of ergodic measures with regular Lyapounov exponents:

\begin{corollary}\label{c.recurrence}Let $\eta$ be an ergodic measure with
$\lambda^c(\eta)>-\infty$. Then, for any $\gamma>0$, there exists $\delta>0$
such that
\begin{equation}\label{e.slow}
\limsup\limits_{n\to\infty}\frac{1}{n}\sum\limits_{j=0}^{n-1} 
-\log \text{ dist}_{\delta} (f^j(p),\SC)\leq\gamma, \quad \text{ for
}\eta-\text{a.e. point } p.
\end{equation}
Here $\text{dist}_{\delta}$ is the $\delta$-truncated distance defined by
\begin{equation}
\text{ dist}_{\delta} (p,\SC)= \left\{ \begin{array}{ll}
1 & \textrm{if $\text{ dist }(p,\SC)\geq\delta$,}\\
\text{dist }(p,\SC) & \textrm{otherwise}
\end{array} \right.
\end{equation}
\end{corollary}

\begin{proof}By the lemma~\ref{l.recurrence}, the function 
$|\log(dist(p,\SC))|$ is $\eta$-integrable. In particular, $\eta (\SC)
= 0$. Using these facts together with the definition of the $\delta$-truncated
distance $dist_{\delta}$, we have that for any $\gamma>0$, there exists
some $\delta>0$ with
$$
\int -\log dist_{\delta} (p,\SC) \,d\eta \leq \gamma.
$$
Because $\eta$ is ergodic, a simple application of Birkhoff's theorem to the
function $-\log(dist_{\delta}(p,\SC))$ completes the proof.
\end{proof}

At this point, we define 

\begin{definition}Let $K$ be the set of ergodic measures $\mu$ whose central 
Lyapounov exponent $\lambda^c (\mu)$ are greater than
$\frac{1}{4}\bigl(c_0 -\frac{\log (d+\ep)}{\log (d-\ep)}\bigr)$.
We define $\SK$
as the set of invariant measures $\mu$ whose ergodic decomposition $(\mu_p)$ 
belongs to $K$ for $\mu$-almost every $p$. The sets $K$ and $\SK$ are
not empty since $\mu_0$ belongs to both of them. See 
theorem~\ref{t.exponent}.
\end{definition}

It is interesting to consider $\SK$ since it contains any measure liable of
satisfying the variational principle:

\begin{lemma}\label{l.K}There exists a constant $\ka_0>0$ such that every 
measure $\eta\notin\SK$ satisfies
\begin{equation}\label{e.K}
h_{\eta} (f) + \int\phi\, d\eta +\ka_0 < 
\sup\limits_{\mu} \left( h_{\mu} (f)+\int\phi\, d\mu\right).
\end{equation}
\end{lemma} 

\begin{proof}The idea is to compare $h_{\eta}(f)$ with $h_{\mu_0}(f)$, where
$\mu_0$ is the unique SRB measure of $f$ (see theorem~\ref{t.SRB}).
Without loss of generality, we can suppose that $\eta$ is ergodic.
In this case, $\eta\notin\SK$ means 
$\lambda^c(\eta)\le\frac{1}{4}\bigl(c_0 -\frac{\log (d+\ep)}{\log(d-\ep)}\bigr)$. 
Ruelle's inequality combined with the
estimate~(\ref{e.conformal}) of proposition~\ref{p.conformal} gives
$$h_{\eta} (f)\leq \log (d+\ep) + \frac{1}{4}\left(c_0 -\frac{\log (d+\ep)}{\log(d-\ep)}\right).$$
Hence, the condition on the potential~(\ref{e.potential}) and
estimate~(\ref{e.1}) imply
$$h_{\eta} (f)+\int\phi\, d\eta < h_{\mu_0} (f) +\int\phi\, d\mu_0 - 
\frac{1}{4}\log (d+\ep).$$
Taking $\ka_0 = \frac{1}{4}\log d < \frac{1}{4}\log(d+\ep)$, for instance,
concludes the proof.
\end{proof}

Before proceeding further in the proof of theorem~\ref{t.A}, we recall the
concept of hyperbolic time for maps with critical points. The setting of this
definition is as follows.

Consider $f:M\to M$ a $C^2$ map which is local diffeomorphism except at a zero
Lebesgue measure set $\SC\subset M$. Assume $f$ \emph{behaves like a power of
the distance} to the critical set $\SC$, i.e., there are constants $B>1$ and
$\ell>0$ such that, for every $p,q\notin\SC$ with
$2\, dist(p,q)<dist(p,\SC)$ and $v\in T_p M$ we have:

\begin{itemize}
\item $\frac{1}{B}\ dist(p,\SC)^{\ell}\leq \frac{\|Df(p)v\|}{\|v\|}\leq
B\ dist(p,\SC)^{\ell}$;
\item $|\log\|Df(p)^{-1}\|-\log\|Df(q)^{-1}\||\leq
B\frac{dist(p,q)}{dist(p,\SC)^{\ell}}$;
\item $|\log|\det Df(p)^{-1}|-\log|\det Df(q)^{-1}|\leq 
B\frac{dist(p,q)}{dist(p,\SC)^{\ell}}$,
\end{itemize}

\begin{remark}\label{r.1}As the reader can easily verify, Viana maps satisfy 
all the previous assumptions, since these conditions are $C^2$-open and
$\widetilde{f}$ satisfies them.
\end{remark}

For the definition of hyperbolic times, we fix $0<b<\min\{1/2,1/(2\be)\}$.
\begin{definition}\label{d.times}
Given $0 < \sigma < 1$ and $\delta > 0$, we say that $n$ is a
$(\sigma,\delta)$-hyperbolic time for $p$ if, for all $1\leq k\leq n$,
\begin{equation*}
\prod\limits_{j=n-k}^{n-1}\|Df(f^j(p))^{-1}\|\leq\sigma^k \quad \text{ and }
dist_{\delta}(f^{n-k}(p),\SC)\geq \sigma^{bk}.
\end{equation*}
\end{definition}

The usefulness of hyperbolic times is the key property:

\begin{proposition}\label{p.expansivity}
Given $\sigma<1$ and $\delta>0$, there exists $\delta_1>0$ such that if $n$ is a
$(\sigma,\delta)$-hyperbolic time of $p$, then there exists a neighbourhood
$V_p$ of $p$ such that:
\begin{itemize}
\item $f^n$ maps $V_p$ diffeomorphically onto the ball $B_{\delta_1}(f^n(p))$;
\item For every $1\leq k\leq n$ and $y,z\in V_p$,
$$\text{dist}(f^{n-k}(y),f^{n-k}(z))\leq 
\sigma^{k/2}\text{dist}(f^{n}(y),f^{n}(z)).$$
\end{itemize}
\end{proposition}

\begin{proof}See Alves, Bonatti and Viana~\cite[p.377]{ABV}.
\end{proof}

Next, we recall the following criterion to prove the existence of infinitely
many hyperbolic times.

\begin{notation}For a fixed $\sigma <1$, let $H(\sigma)$ be the set of points
$p\in M$ with the following two properties:

$$\limsup\limits_{n\to\infty}\frac{1}{n}\sum\limits_{j=0}^{n-1}
\log\|Df(f^j(p))^{-1}\|\leq 3\log\sigma <0 $$
and, for any $\ga>0$ there is some $\delta>0$ satisfying 
$$\limsup\limits_{n\to\infty}\frac{1}{n}\sum\limits_{j=0}^{n-1} 
-\log\bigl(dist_{\delta}(f^j(p),\SC)\bigr)\leq\gamma.$$
\end{notation}

\begin{proposition}\label{p.infinitely}
Given $\sigma<1$, there exist $\nu>0$ and $\delta>0$ depending only on
$\sigma$ and $f$ such that, given any $p\in H(\sigma)$ and $N\geq 1$
sufficiently large, there are 
$(\sigma,\delta)$-hyperbolic times $1\leq n_1<\dots<n_l\leq N$ for $p$ with
$l\geq \nu N$.
\end{proposition}

\begin{proof}See Alves, Bonatti and Viana~\cite[p.379]{ABV}
\end{proof}

Coming back to the proof of theorem~\ref{t.A}, we introduce the following 
definition:

\begin{definition}\label{d.expanding}$\mu$ is called a 
$\sigma^{-1}$-\emph{expanding measure} if $p\in H(\sigma)$ for $\mu$-almost
every point $p$. 
\end{definition}

We now prove that any measure $\mu\in\SK$ is expanding:

\begin{lemma}\label{l.expanding}
Any $\mu\in\SK$ is a $\sigma^{-1}$-expanding with $\sigma = \exp
(-\frac{1}{12}\zeta)$, where $\zeta:=c_0-\frac{\log (d+\ep)}{\log (d-\ep)}>0$.  
\end{lemma}

\begin{proof}Without loss of generality we can assume $\mu\in\SK$ ergodic. This
implies that 
$$\lambda^c(p) = \lambda^c(\mu) = \int\log\|Df|_{E^c}\| d\mu\geq
\frac{1}{4}\left(c_0-\frac{\log (d+\ep)}{\log (d-\ep)}\right)$$
for $\mu$-a.e. $p$. Because $\mu$ is ergodic and the central direction $E^c$ is
one-dimensional, this is equivalent to 
$$
\lim_{n\to\infty}\frac{1}{n}\sum_{j=0}^{n-1}\log\|Df(f^j(p))^{-1}\|\leq 
-3 \frac{1}{12}\left(c_0-\frac{\log (d+\ep)}{\log (d-\ep)}\right) = 3\log\sigma < 0.$$
This fact and the corollary~\ref{c.recurrence} gave us the assumptions in 
the definition of $H(\sigma)$ at $\mu$-almost every point.
\end{proof}

An easy consequence of the expanding features of a given measure is the
existence of generating partitions, which turns out to be a relevant property
when studying entropy and equilibrium states:

\begin{lemma}\label{l.partition}
Given $\sigma<1$, there exists $\delta_1>0$ such that any 
partition $\SP$ with diameter less than $\delta_1$ is a generating partition for
any $\sigma^{-1}$-expanding measure $\mu$.
\end{lemma} 

\begin{proof}
Define 
$$
A_{\ep}(p):=\{y\in M: dist(f^j(y),f^j(p))\leq\ep \quad \text{for every }n\ge
0\}.
$$
First we prove that $A_{\ep}(p)=\{p\}$ for $\mu$-a.e. $p$
and then we show how this can be used to finish the proof.
\begin{claim}\label{claim}
For any $\sigma^{-1}$-expanding measure $\mu$, there exists $\delta_1>0$ such
that for any $\ep<\delta_1$ and $\mu$-a.e. $p$,
$$A_{\ep}(p)=\{p\}.$$
\end{claim}

\begin{proof}[Proof of the claim]
Proposition~\ref{p.infinitely} guarantees that $\mu$-a.e. $p$ has infinitely
many hyperbolic times $n_i(p)$, since $\mu\in\SK$. Hence, applying the
proposition~\ref{p.expansivity} we conclude that there exists some $\delta_1>0$
such that if $z\in A_{\ep}(p)$ with $\ep<\delta_1$ then for any $n_i$ we have 
$$
dist(z,p)\leq \sigma^{n_i/2}dist(f^{n_i}(z),f^{n_i}(p))\leq \sigma^{n_i/2}\ep.
$$
Since $n_i\to\infty$ as $i\to\infty$, we deduce that $z=p$.
\end{proof}

It is now a more or less standard matter to show that the previous claim
implies that any finite partition $\SP=\{P_1,\dots,P_l\}$ with
$\text{diam}(\SP)<\delta_1$ is a generating partition.
Indeed, given any measurable set $A$ and given $\delta>0$, consider $K_1\subset
A$ and $K_2\subset M-A$ compact sets with $\mu(A\setminus K_1)\leq\delta$ and
$\mu(A^c\setminus K_2)\le\delta$. 
Let $r:=dist(K_1,K_2)>0$. Claim~\ref{claim} shows that, for $n$ large
enough, $\text{diam}\SP^{(n)}(p)\leq r/2$ for every $p$ in a set of 
$\mu$-measure greater than $1-\delta$, where
$$
\SP^{(n)}(p):=
\{C^{(n)}:=P_{i_1}\cap\dots f^{-n+1}(P_{i_n})\mbox{ with }p\in C^{(n)}\}.
$$
Consider the sets $C_1^{(n)},\dots,C_m^{(n)}\in\SP^{(n)}$ that intersect $K_1$.
Then, it is not difficult to see that
$$\mu(\cup C_i^{(n)}\Delta A)\leq 3\delta.$$
Since $A$ is an arbitrary measurable set and $\delta>0$ is an arbitrary positive
real number, this proves that $\SP$ is a generating partition.
\end{proof}

Finally, we conclude the proof of theorem~\ref{t.A}:

\begin{proof}[Proof of theorem~\ref{t.A}]We take a sequence $\mu_k$ of invariant
measures such that $h_{\mu_k}(f)+\int\phi
d\mu_k\to\sup_{\eta}(h_{\eta}(f)+\int\phi d\eta)$. when $k\to\infty$. 
Taking a subsequence if necessary, we also assume $\mu_k\to\mu$ as $k\to\infty$.
By the lemma~\ref{l.K} we may assume that $\mu_k\in\SK$.

We are going to prove that $\mu$ is an equilibrium measure and every equilibrium
measure belongs to $\SK$. Clearly, these claims finish the proof.

Fix a finite partition $\SP$ with diameter less than $\delta_1>0$, the
constant of lemma~\ref{l.partition} (with $\sigma$ being the constant of 
lemma~\ref{l.expanding}) such that $\mu (\p P) = 0$ for any $P\in\SP$. Observe
that Kolmogorov-Sinai's theorem implies $h_{\mu_k}(f)=h_{\mu_k}(f,\SP)$ by the
lemmas~\ref{l.partition} and~\ref{l.expanding}, since $\mu_k\in\SK$. Because the
function $\nu\to h_{\nu}(f,\SP)$ is upper-semicontinuous at every measure $\nu$
with $\nu(\p\SP)=0$, we get
\begin{eqnarray*}
\sup_{\eta}(h_{\eta}(f)+\int\phi d\eta) &=& \limsup\limits_k 
(h_{\mu_k}(f)+\int\phi d\mu_k) \\ &=&  
\limsup\limits_{k} (h_{\mu_k}(f,\SP)+\int\phi
d\mu_k)\\ &=& h_{\mu}(f,\SP)+\int\phi d\mu \\ 
&\leq& h_{\mu}(f)+\int\phi d\mu.
\end{eqnarray*}
In other words, $\mu$ is an equilibrium state. Finally, the fact that 
every equilibrium state $\eta$ belongs to $\SK$ follows directly from the
lemma~\ref{l.K}. This concludes the proof.
\end{proof}

\section{Remarks}\label{s.remarks}

Let us make few comments about the hypothesis $d\geq 16$ in the definition of
the Viana maps, about the generalization of the theorem~\ref{t.A} in
the context of random perturbations and about some natural questions associated
to its statement.

The hypothesis $d\geq 16$ in the definition of the maps $\widetilde{f}$ was
used by Viana to ensure that $C^3$ small perturbations have two positive
Lyapounov exponents. However, Buzzi, Sester and Tsujii~\cite{BST} noticed that 
if one is willing to consider $C^{\infty}$ small perturbations (instead of $C^3$), then
the theorems of Viana and Alves (see theorems~\ref{t.exponent},~\ref{t.SRB})
still hold. Also, it is proved in~\cite{BST} that the vertical expansion on the
dynamical strip $I$ of the derivative of $f$ is uniformly less than $2$. In
particular, our comments on the tangent bundle dynamics of $f$ are true even if
$2\leq d\leq 16$. Summarizing,
\begin{center}
\emph{Theorem~\ref{t.A} also holds for $C^{\infty}$ small perturbations of
  Viana maps with} $d\geq 2$
\end{center}
 (instead of stronger requirement $d\geq 16$).

Also, it is not difficult to check that the arguments of this note apply to the
context of small random perturbations of Viana maps. Indeed, the techniques
of Alves and Ara\'ujo~\cite{AA} allow one to consider hyperbolic times,
expanding measures, slow recurrence and non-uniform expansion in the
non-deterministic setting
\footnote{Although the statements in Alves and Ara\'ujo~\cite{AA} are
for Viana maps (and its $C^3$ perturbations) with $d\geq 16$,  the results of
Buzzi, Sester and Tsujii~\cite{BST} ensure their proof can be carried out for
$d\geq 2$ (at the cost of considering only $C^{\infty}$ small perturbations).}
, and prove the analogues (in the random context) of Viana's result (see
theorem~\ref{t.exponent}) and Alves' result (see theorem~\ref{t.SRB}) about
the existence of a unique physical measure with all Lyapounov exponents
larger that some positive number.
We obtain a version of theorem~\ref{t.A} for small random perturbations of
Viana maps by considering the modified definition of the set of measures $\SK$ 
(as defined in Arbieto, Matheus, Oliveira~\cite{AMO}). 

Finally, a natural question motivated by the existence of equilibrium measures
is the uniqueness and ergodic propoerties, including decay of correlations,
for the equilibrium state obtained in theorem~\ref{t.A}.
Since the approach of Young towers is not easy in the context of general
equilibrium measures\footnote{This occurs essentially because 
the known constructions of Young towers depends on the Lebesgue measure as a
reference measure, but the general equilibrium measures are singular with
respect to Lebesgue, even in the uniformly hyperbolic case.} and Ruelle's
thermodynamical formalism via the spectral properties of the transfer operator
is delicate in the presence of critical points, this question presents an
interesting problem for the extension of the theory of equilibrium states
beyond uniform hyperbolicity. For results in this direction in one dimension,
see \cite{Pesin-Senti1}, \cite{Pesin-Senti2}. 

\textbf{Acknowledgments.} The authors are thankful to IMPA and its staff for
the fine reseach ambient. A.A. and C.M. would like to acknowledge Krerley
Oliveira for pointing out a gap in an early version of this note. Also,
A.A. and C.M. are indebted to Viviane Baladi for the invitation for a
post-doctoral visit at the Institut de Math\'ematiques de Jussieu where this 
project started and for the trimester ``Time at Work'' at the Institute Henri 
Poincar\'e where this paper was completed. A.A. and C.M. were partially
supported by the Conv\^enio Brasil-Fran\c ca em Matem\'atica. 
S. S. was supported by a grant by the Swiss National Science Foundation.
    
\bibliographystyle{alpha}
\bibliography{bib}

\vspace{2.5cm}

\noindent    \textbf{Alexander Arbieto} ( alexande{\@@}impa.br )\\
             \textbf{Carlos Matheus} ( matheus{\@@}impa.br )\\
	     \textbf{Samuel Senti} ( senti{\@@}impa.br )\\
             IMPA, Est. D. Castorina 110, Jardim Bot\^anico, 22460-320 \\
             Rio de Janeiro, RJ, Brazil

\end{document}